\documentclass{amsart}

\begin{document}
\title{Path  algebras, wave-particle duality, and quantization of phase space}
\author{Murray Gerstenhaber}
\address{Department of Mathematics\\  University of Pennsylvania\\
  Philadelphia, PA 19104-6395, USA}
\email{mgersten@math.upenn.edu}
\thanks{The author thanks V. Dolgushev, J. D. Stasheff, and  D. Sternheimer for many helpful suggestions.}
\subjclass[2010]{16S80 (primary); 53B50, 53D55 (secondary)}

\begin{abstract}\noindent 
Semigroup algebras admit certain `coherent' deformations which, in the special case of a path algebra, may associate a periodic function to an evolving path; for a particle moving freely on a straight line after an initial impulse, the wave length is that hypothesized by de Broglie's wave-particle duality. This theory leads to a model of  ``physical'' phase space of which mathematical phase space, the cotangent bundle of configuration space, is a projection.  This space is singular, quantized at the Planck level, its structure implies the existence of spin, and the spread of a packet can be described as a random walk. The wavelength associated to a particle moving in this space need not be constant and its phase  can change discontinuously.\end{abstract}

\keywords{cohomology, deformation, homotopy, semigroups}
\dedicatory{In memory of my mother, a seamstress}

\maketitle
\newtheorem{theorem}{Theorem}
\newtheorem{corollary}{Corollary}
\newtheorem{lemma}{Lemma}
\newtheorem{defn}{Definition}
\renewcommand{\abstractname}{}
\newcommand{\C}{\ensuremath{\mathbb{C}}}
\newcommand{\R}{\ensuremath{\mathbb{R}}}
\newcommand{\Z}{\ensuremath{\mathbb{Z}}}
\newcommand{\cP}{\ensuremath{\mathcal{P}}}
\newcommand{\cM}{\ensuremath{\mathcal{M}}}       
\newcommand{\I}{\ensuremath{\mathcal{I}}}
\newcommand{\bM}{\ensuremath{\mathbf{M}}}
\newcommand{\bS}{\ensuremath{\mathbf{S}}}
\newcommand{\cX}{\ensuremath{\mathcal{X}}}
\newcommand{\cA}{\ensuremath{\mathcal{A}}}
\newcommand{\bA}{\ensuremath{\mathbf{A}}}
\newcommand{\bR}{\ensuremath{\mathbf{R}}}
\newcommand{\cB}{\ensuremath{\mathcal{B}}}
\newcommand{\cS}{\ensuremath{\mathcal{S}}}
\newcommand{\pr}{\ensuremath{\preceq}}
\newcommand{\op}{\ensuremath{\mathrm{op}}}
\newcommand{\bfk}{\ensuremath{\mathbf{k}}}
\newcommand{\g}{\gamma}
\newcommand{\om}{\omega}
\newcommand{\Om}{\Omega}
\newcommand{\G}{\Gamma}
\newcommand{\z}{\zeta}
\newcommand{\Rp}{{R}^+}
\newcommand{\Aut}{\operatorname{Aut}}
\newcommand{\dr}{\operatorname{dR}}
\vspace{-7mm}
\date{}

Semigroup algebras admit certain `coherent' deformations which, in the special case of a path algebra, may associate a periodic function to an evolving path.  This stems from the fact that if $A$ is a semigroup algebra, then the Hochschild complex $C^*(A,A)$ contains a distinguished subcomplex consisting of its coherent cochains, defined below.
 In phase space
there is a canonical deformation of the path algebra induced by the symplectic form. For a particle moving freely on a straight line after an initial impulse, the periodic function associated to it by the deformation has wave length that hypothesized by de Broglie's wave-particle duality.

Algebraic deformation theory, introduced in the author's paper \cite{G:Def1}, has found an increasing number of important applications, amongst them by Bayen et al,  \cite{BFFLS},  to the creation of an autonomous phase space approach to quantization  (expressing Moshe Flato's vision that quantization is deformation\footnote{Groenewold, \cite{Groenewold}, and
Moyal, \cite{Moyal}, implicitly used special instances of  the universal deformation formula associated to the two-dimensional Abelian Lie algebra, \cite[p. 13]{G:Def3}, \cite{GiaquintoZhang:Bialgebra}. For origins and an overview of deformation quantization, cf
\cite{DitoSternheimer:Genesis}, \cite{Sternheimer:20YearsAfter}; for a comprehensive treatise, cf \cite{Waldmann:Poisson-Geometrie}. }), by Kontsevich, \cite{Kontsevich:Poisson}, to the deformation of Poisson manifolds, by Rieffel to strict quantization of $C^*$ algebras, cf \cite{Rieffel:actions}, and more recently by Romero-Ayala et al, \cite{Romero-Ayala}, to the study of black holes and general relativity.

Quantization has generally been viewed as deformation of a commutative multiplication to one which is no longer commutative. This paper extends the deformation approach to quantization by applying it to cases where multiplication is initially non-commutative.  Planck's original quantum of action, $h \approx 6.626070\times10^{-34}$ Joule$\cdot$seconds is a constant of nature but is treated as a variable in deformation theory. (It is frequently convenient to absorb $2\pi$ into $h$, leading to the notation  $\hbar =h/2\pi \approx 1.0545718\times10^{-34} \text{ J}\cdot$s; more recently, $\hbar$ has also become a common notation for a deformation parameter.) 
Classical deformation quantization is the theory of \cite{G:Def1} applied to products of functions in phase space, most particularly, of position and momentum. When $h = 0$, these commute, but cease to commute when $h > 0$. The `coherent' deformation theory of semigroup algebras introduced here is meaningful, in particular, for  products (concatenations) of paths, which is naturally non-commutative.  
When $h = 0$ 
one has the usual concatenation.  With $h > 0$, however,  
the product may be multiplied by a periodic factor, in effect quantizing the path algebra and attaching a wave to a moving particle. A particle set in motion by an impulse (this being, in effect, the concatenation of two orthogonal paths in phase space), will have then have an attached wave length of $h/p$, where $p$ is the momentum of the particle, as predicted by de Broglie and first observed in the Davisson-Germer experiment.

Section \ref{twist} presents basic results on the cohomology of semigroups.  While  important in themselves, for the later sections one will need only certain definitions, namely that of the subcomplex of coherent cochains inside the Hochschild cochain complex of a semigroup algebra, as well as that of coherent deformations  and  twisting functions (multiplicative 2-cocycles) of a semigroup.  In \S \ref{path} Moore paths are reviewed. These produce a strictly associative multiplication on paths, the usual multiplication of paths being only homotopy associative. 
 
Section \ref{taille} introduces the concept of \emph{taille}:  A $k$-cocycle $\omega$ of a space $X$ with coefficients in an Abelian group $G$ determines a morphism of its $k$th homotopy group $\pi_k$ into $G$ which depends only on the class of $\omega$ in $H^k(X, G$).  The image group, $\omega(\pi_k)$, will be called the \emph{taille group of $X$ relative to $\omega$}.  The main case of interest will be that where  $X$ is a manifold $\cM$, in which case $G = \R$, and where $k = 2$. If the taille group of $\cM$ relative to a closed de Rham 2-form $\omega$ is a discrete subgroup of $\R$,
 then the \emph{taille of $\cM$ relative to $\omega$}, briefly the taille of $\omega$, will  be understood to be  zero if the group is trivial, its positive generator if the group is non-trivial, and $\infty$ otherwise. When $\cM$ is symplectic, by its \emph{taille} we will mean  its taille  relative to its defining 2-form. If $\cM$ is Riemannian,
  then from a closed de Rham 2-form $\omega$ with finite taille $\tau$ one can construct  an additive \emph{tailleur 2-cocycle}  $\tilde\omega$ of the path semigroup $\bS\cM$ of $\cM$. The coefficients of $\tilde\omega$ will be in $\R$ if $\tau = 0$ and in $\R/\tau$ if $\tau > 0$. If $\tau = 0$ then $\omega$  induces a real one-parameter family of deformations of the path algebra of $\cM$. The deformation is then defined by a \emph{twisting function} with positive real values, obtained by exponentiating  the tailleur  2-cocycle.  If $\tau > 0$  then the twisting function takes its values in the circle group; it becomes periodic as a function of time for a free particle moving in flat phase space. The twisting function  defines  a single new associative  multiplication on the path  algebra which is, in effect,  a quantization of that algebra.  An additive 2-cocycle with values in $\R$ can be reduced modulo any $\tau > 0$ to give one  with values in $\R/\tau$ but the deformation produced by this reduced cocycle behaves  differently from that given by the original, something important for wave-particle duality.

Quantization of phase space is introduced in  \S \ref{dB}. 
Mathematical phase space, the cotangent bundle of configuration space, always has taille equal to zero relative to its canonical symplectic 2-form   $\omega =\sum dp_i\wedge dq^i$ where the $q^i$ are configuration coordinates and the  $p_i$ their conjugate momenta, while physical indications require that the taille of phase space be $h$. Since the taille of the cotangent bundle is zero,
  the tailleur cocycle can, as we will see, be reduced to having values in $\R/h$, which then gives the  wave length  predicted by de Broglie. (For this, it will be necessary in \S \ref{dB} to use the Maslov index, \cite{Arnol'd:Maslov}.) This reduction evades, however,  the need to give phase space a  structure with which it has an intrinsic taille of  $h$ rather than zero. To resolve this paradox, we propose,
   for the simplest case of a particle moving on a line, a model of a `physical' phase space of which mathematical phase space is (except for a set of measure zero) a 2 to 1 projection. This physical phase space is quantized, being constructed at the Planck scale of discrete but connected building blocks. As a result it has singularities but also introduces spin, suggesting that any model of physical phase space which intrinsically gives it a positive taille must  introduce spin.  The section concludes with some speculations about the spread of a packet and about the nature of wave-particle duality in a possibly curved phase space.

\section{Semigroup algebras, coherent cochains, and twists}\label{twist}
\noindent A semigroup $\bS$ is a set with an associative multiplication. If it has a unit element, which is then unique and denoted by $I$ or 1, it is called a monoid. When $\bS$ has an absorbing element,  also necessarily unique and  denoted by 0,  it is a ``semigroup with 0'' (or  monoid with 0 if it has both).  Together with any commutative unital ring $\bfk$, one can form the \emph{semigroup algebra} $\bfk\bS$, which is the free module over $\bfk$ generated by the elements of  $\bS$ with multiplication determined by that in $\bS$.  

Polynomial algebras, commutative and non-commutative, are semigroup algebras with underlying semigroup $\bS$ consisting of their monomials. Certain quotients of these, e.g., truncated polynomial algebras,  are again semigroup algebras. Matrix algebras form another class of examples, where $\bS$ consists of the matrix units $e_{ij}$ together with the zero element. If $\mathcal I = \{i,j,...\}$ is a poset with partial order denoted by $\preceq$, then its poset algebra $\bfk\mathcal I$  (the free module generated over $\bfk$ by formal  elements $e_{ij},\, i\preceq j$ with multiplication given by $e_{ij}e_{kl} = \delta_{jk}e_{il}$)  is a semigroup algebra;  the underlying semigroup $\bS$ is again the set of all $e_{ij}$ together with zero. When $\mathcal I$ is finite,
 one can introduce a linear order compatible with the partial order; $\bfk\mathcal I$ can then be viewed as an algebra of upper triangular matrices and the $e_{ij}$ again as matrix units. An algebra over a field having a basis in which the product of any two elements is either another basis element or zero is said to have a multiplicative basis. Any such algebra is clearly also a semigroup algebra. Representation finite algebras, i.e., ones admitting only a finite number of isomorphism classes of indecomposable modules, have such bases, \cite{FiniteRep}. The semigroup algebras of principal interest here will be path algebras, \S \ref{path}, with multiplicative bases consisting of their paths. They will generally have no unit element, although one can be adjoined if necessary. Amongst the common path algebras are quiver algebras. On an arbitrary topological space, the Moore paths form a semigroup; the corresponding semigroup algebra with coefficients in the  complex numbers will be called \emph{the} path algebra of the space. 
 
What all the above have in common is that every semigroup algebra $A$ contains a distinguished subcomplex of its Hochschild cochain complex $C^*(A, A)$ consisting of its coherent cochains:
\begin{defn}\label{coherent}
Let $A = \bfk\bS$ be a semigroup algebra and $C^n(A,A)$ be its module of $n$-cochains with coefficients in itself. 
An $n$-cochain $f\in C^n(A,A),\, n > 0$ will be called \emph{coherent}  if $f(a_1, \dots, a_n)$ is always a multiple of $a_1a_2\cdots a_n$ by some element of $\bfk$. With the convention that  when $a_1a_2\cdots a_n = 0$  this multiple is taken to be zero, we can then write $f(a_1, \dots , a_n ) = F(a_1, \dots , a_n)a_1a_2\cdots a_n$ for some uniquely defined function $F$ associated to $f$.
\end{defn}

In this section we will state without proof some of the more readily demonstrated propositions which are not needed in the rest of this paper.

The module of coherent $n$-cochains of $A$ will be denoted by $CC^n(A,A)$; we set $CC^0(A,A) = 0.$ The Hochschild coboundary of a coherent cochain is clearly again coherent.  The coherent cochains therefore form a subcomplex $CC^*(A,A)$ of the full Hochschild complex $C^*(A,A)$, so one can define the \emph{coherent cohomology} of $A$,  denoted  $CH^*(A,A)$.   Cup products of coherent cochains are coherent.  Moreover, if $f^m, g^n$ are coherent cochains of dimensions $m,n$ respectively, then their composition product $f^m\circ_i g^n \in CC^{m+n-1}$, given by
$$f^m\circ_i g^n(a_1,\dots,a_{m+n-1}) = f^m(a_1,\dots,a_{i-1}, g^n(a_i,\dots, a_{i+n-1}), a_{i+n}\dots,a_{m+n-1})$$ 
is coherent, and from these compositions $CC^*(A,A)$ acquires a graded pre-Lie structure. As the coherent cochain complex $CC^*(A,A)$ is thus closed under all the operations used in \cite{G:Cohomology}, we have the following.

\begin{theorem}
If $A= \bfk\bS$ is a semigroup algebra then its coherent cohomology,
 $CH^*(A,A)$, is a Gerstenhaber algebra. \quad $\square$
\end{theorem}
 
In the particular case when $A$ is a commutative algebra over a field $\bfk$ of characteristic zero and $M$ a symmetric $A$-bimodule, one has a Hodge-type decomposition, \cite{GS:Hodge}, of the Hochschild cochain complex $C^*(A, M)$ into an infinite direct sum of subcomplexes defined by idempotents in the group ring of the symmetric group. These idempotents, originally called  either the `Gerstenhaber-Schack' (GS) or `homological' idempotents, were later designated `Eulerian' idempotents; we will call them here the GS idempotents. When $\bfk$ has characteristic $p>0$, 
the decomposition of $C^*(A,M)$ collapses into a direct sum of $p-1$ subcomplexes, \cite{Loday:Euler}. These decompositions induce corresponding ones of the cohomology, $H^*(A,M)$.  In a commutative semigroup algebra, the GS idempotents carry coherent cochains into coherent cochains. These decompositions therefore also restrict to $CC^*(A,A)$, so we have the following.

\begin{theorem}
Let $\bS$ be a commutative semigroup, and set $A = \bfk\bS$, where $\bfk$ is a field.  If the characteristic of $\bfk$ is zero then $CC^*(A,A)$  has a full Hodge-type decomposition; if it is $p >0$ then $CC^*(A, A)$, decomposes into a direct sum of $p-1$ subcomplexes. The decompositions descend to $CH^*(A,A).$ \quad $\square$
\end{theorem} 

Let $\mathcal {I}$ be a (possibly infinite) poset, $N\mathcal{I}$ be its nerve, and $A = \bfk\mathcal{I}$ be its poset algebra with coefficients in the commutative unital ring $\bfk$; $A$ is quasi self-dual in the sense of \cite{G:Self-dual}, and in the category of such algebras, $H^*(A,A)$ is a contravariant functor.
There is then an isomorphism,  natural in $\mathcal{I}$, $H^*(N\mathcal{I}, \bfk) \to H^*(\bfk\mathcal{I},\bfk\mathcal{I})$, cf \cite{GS:SC=HC} . 
If $K$ is a simplicial complex then its barycentric subdivision $K'$ has the same cohomology as $K$ and is also a poset, so its poset algebra has the same cohomology as $K$. The following is essentially a restatement of the main result of \cite{GS:SC=HC}.

\begin{theorem}\label{poset}
If $A$ is a poset algebra, then the inclusion of $CC^*(A,A)$ into $C^*(A,A)$ induces an isomorphism of cohomology $CH^*(A,A) \cong H^*(A,A)$. $\square$
\end{theorem}
The Hochschild-Kostant-Rosenberg theorem implies that if $A = \bfk[x_1,\dots, x_N]$ is a polynomial ring over a field $\bfk$ of characteristic zero,  then $H^*(A,A)$ can be identified with the exterior algebra on the derivations of $A$.   
\begin{theorem}Let $A = \bfk[x_1,\dots,x_N]$  be a polynomial ring in $N$ variables over a field $\bfk$ of characteristic zero. Then $CH^*(A,A)$ can be identified with the exterior algebra on the space spanned by the coherent derivations $x_1\partial_{x_1}, \dots, x_N\partial_{x_N}$. $\square$
\end{theorem}
The Hochschild-Kostant-Rosenberg theorem and the foregoing  must both be modified when $\bfk$ has characteristic $p > 0$ since the isomorphism $H^*(A,A)\cong \bigwedge \operatorname{Der} A$ is then only one of vector spaces, but we will not need that here.

A \emph{coherent deformation} of a semigroup algebra $\bfk\bS$ is one in which all the cochains involved are coherent. For  example, if $A = \C[x,y]$, then $H^*(A,A)$ is spanned by the single cohomology class of  $\partial_x\wedge \partial_y = (1/2)(\partial_x\smile \partial_y - \partial_y\smile \partial_x)$, which is not coherent (but is the infinitesimal of the jump deformation to the first Weyl algebra), while $CH^*(A,A)$ is spanned by the single cohomology class of $x\partial_x\wedge y\partial_y.$ The star product arising from this cocycle deforms $\bfk[x,y]$ to the quantum plane, for writing $q$ for $e^{\hbar}$, where $\hbar$ is the deformation parameter, one has
$x*y = e^{q/2}xy, \,y*x = e^{-q/2}xy$, so $x*y = e^qy*x$.  (The last equation completely defines the structure of the deformed algebra but not the deformation. One would obtain the same algebra, up to isomorphism, by using any cocycle cohomologous to $x\partial_x\wedge y\partial_y$, e.g., $x\partial_x\smile y\partial_y$,  
but the physical implications of the deformation will be different, cf, e.g., \cite{G:Least uncertainty}.)

The function $F$ which appears in Definition \ref{coherent} is just an $n$-cochain of the semigroup $\bS$ with coefficients in $\bfk$. 
If $\bfk$ is any commutative unital coefficient ring, then one can define the cohomology $H^*(\bfk\bS, M)$ with coefficients in any $\bfk\bS$-bimodule $M$, but if $\bS$ has a zero element,
 then $H^n(\bfk\bS, M)$ will vanish for all $n > 0$.  In this case, therefore, one must restrict cochains $F^n \in C^n(\bS, M)$ by letting them be  defined only on the set $\bS^{(n)}$ of  those $n$-tuples $(a_1, \dots, a_n) \in \bS^n$ with $a_1a_2\cdots a_n \ne 0$, cf \cite{Novikov:Semigroup}. One can then take $M$ to be an arbitrary Abelian group $G$, with the operation of all non-zero elements of $\bS$ being the identity. 
Taking $G$ to be additive, for all $n > 0$ define the \emph{additive} $n$-cochains $C^n(\bS, G)$ of $\bS$ with coefficients in $G$ to be formal finite sums with integer coefficients of mappings $F\!:\! \bS^{(n)} \to G$. The coboundary operator $\delta: C^n \to C^{n+1}$ is given by 
\begin{multline}\label{coboundary}
\delta F(a_1,\dots,a_{n+1}) = F(a_2,\dots,a_{n+1})\, +\\ \sum_{i=1}^{n}(-1)^i F(a_1,\dots,a_{i-1}, a_ia_{i+1},a_{i+2}, \dots,a_{n+1}) \\+ (-1)^{n+1}F(a_1,\dots,a_n),
\end{multline}
whenever $a_1a_2\cdots a_n \ne 0$. Then  $\delta\delta F =0$ whenever it is defined on some set of arguments, so setting $C^0 = 0$, cohomology groups can be defined as usual.
 Taking $G = \bfk$, one can identify the complexes $C^*(\bS, \bfk)$ and $CC^*(\bfk\bS,\bfk\bS)$, from which it follows that $H^*(\bS, \bfk)$ inherits the cohomology of the latter. One therefore has the following.

\begin{theorem}
  The cohomology of a semigroup with coefficients in a commutative ring has the structure of a Gerstenhaber algebra. $\square$
\end{theorem}
The cup product and the composition operations in $C^*(\bS, \bfk)$ can be deduced immediately from those of $CC^*(\bfk\bS,\bfk\bS)$: If $F^m$ and
$G^n$ are, respectively, $m$ and $n$ cochains then   
$$
(F^m\smile G^n)(a_1,\dots, a_{m+n}) = F^m(a_1,\dots,a_m) G^n(a_{m+1},\dots,a_{m+n}),$$ 
\vspace{-7mm}
\begin{multline*}
F^m\circ_i G^n(a_1,\dots,a_{m+n-1}) = \\F^m(a_1,\dots,a_{i-1}, a_ia_{i+1}\cdots a_{i+n-1}, a_{i+n},\dots,a_{m+n-1})G^n(a_i,\dots, a_{i+n-1}).
\end{multline*}
In the latter, $G$ simply multiplies $F$, since its values are in $\bfk$.

When the coefficient group $G$ is Abelian but multiplicative, the coboundary formula can be rewritten in multiplicative form and cochains will then be called \emph{multiplicative}. For a multiplicative 1-cocycle $g$ one then has  $g(ab) = g(a)g(b)$ whenever $ab \ne 0$. This is a limited kind of morphism of $\bS$ into $G$ since we clearly can not have  $g(ab) = g(a)g(b)$ when $ab = 0$. For a multiplicative 2-cocycle $f$ one has $f(a,b)f(ab,c) = f(b,c)f(a,bc)$ whenever $abc \ne 0$. 

\begin{defn} A  \emph{twist} of a semigroup algebra $\bfk\bS$ is a multiplicative 2-cocycle of $\bS$ with coefficients in $\bfk^{\times}$.
\end{defn}
A twist $f$ defines a new associative multiplication on  the semigroup $\bfk^{\times}  \times \bS$, and hence on $\bfk\bS$: For  $a, b \in \bS$, set (in simplified notation) $a\ast b = f(a,b)ab$,  whenever $ab \ne 0$, else set $a\ast b = 0$; this is  then extended bilinearly to all of $\bfk\bS$.   The twist $f$ will be called \emph{trivial} if it is a coboundary, i.e.,
 if for some 1-cochain $g$ with values in $\bfk^{\times}$, one has $f(ab) = g(a)g(b)$ whenever $ab \ne 0$.  The mapping of $\bfk\bS$ to itself sending $a \in \bS$ to $g(a)a$ is then a \emph{coherent  isomorphism} of $\bfk\bS$ with the $\ast$ multiplication to $\bfk\bS$ with its original multiplication. A non-trivial twist might produce an algebra isomorphic to  $\bfk\bS$, but not by a coherent isomorphism. Twists constructed from higher dimensional cocycles lead to Stasheff $A_{\infty}$ algebras, cf \cite{Stasheff:Parallel}.

Suppose now that  $\bfk$  is  $\R$ or $\C$ and that we have an additive cocycle $F$ of $\bS$. Then $f_{\hbar}(a,b) = \exp (\hbar F(a,b))$ is a coherent twist for all $\hbar \in \bfk$,  and is embedded in a one-parameter  family of  coherent deformations. Since $f_{\hbar} f_{\hbar'} = f_{\hbar+\hbar'}$, setting $\hbar' = -\hbar$ we can deform back to the original algebra. A deformation constructed in this way can therefore never be a jump deformation, for the infinitesimal of a jump deformation becomes a coboundary in the deformed algebra;  it is not possible to deform back to the original.

Pointwise products of multiplicative cocycles of the same dimension are again cocycles.  In particular,
if $F$ and $F'$ are twisting cocycles, then so is $FF'$ defined by $FF'(a,b) =F(a,b)F'(a,b)$, so one can take integral powers of a twist. When $\bfk = \R$ and the values of the twist are all positive, we can take all real powers, so in this case the twist can be embedded in a coherent deformation. For a nowhere vanishing twist with values in $\C$, 
this is not generally possible; such a twist is the product of a deformation and of a twist whose absolute value is unity. Taking the argument of the latter gives an additive 2-cocycle with coefficients in $\R/2\pi$. 

Suppose that we have a semigroup algebra $\R\bS$ and that we have an additive 2-cocycle $F$ of $\bS$ with coefficients in $\R/\tau$ for some modulus $\tau > 0$.  Extending coefficients to $\C$, we can then define a twist $f$ of $\C\bS$ by setting 
\begin{equation}\label{f}
f(a,b) = \exp((2\pi i/\tau) F(a,b)), \quad  a,b \in \bS.
\end{equation}
Note that if $F$ is an additive cocycle with coefficients in $\R$, then it certainly also defines a cocycle with coefficients in $\R/\tau$, but the twist which the latter defines is quite different. This will make it possible,
in the final section, to construct a deformation of the path algebra of phase space using a modulus $\tau$ equal to    $h$. The canonical symplectic form, multiplied by the Maslov index, \cite{Arnol'd:Maslov},
 will then induce a new twist of its path algebra.  
 
A twist $f$ associated with a positive $\tau$ can no longer be viewed as a classical deformation (even though it may have arisen from an infinitesimal), for in general, a twist of an algebra $A$ need not be part of any parameterized family. It is then not even a `perturbation', which  generally denotes  a parameterized family  which, for some value of the parameter, specializes to $A$.  While the concepts of deformation and perturbation  coincide for finite-dimensional algebras, they do not for infinite dimensional ones. For example, the family of $q$-Weyl algebras  $W_q =\C{\{x,y\}}/(xy-qyx -1)$  specializes to the usual Weyl algebra $W$ at $q = 1$, but $W$ has no cohomology in positive dimensions and in particular is absolutely rigid.  The family of $q$-Weyl algebras is therefore a perturbation but not a deformation of $W$.   For a family of finite dimensional algebras, the dimensions of its cohomology spaces are upper semicontinuous functions of the parameter. This need not be the case when one has a perturbation which is not a deformation, for in that case, some `fragile' cohomology classes may be lost at a given specialization. The cohomology of the specialized algebra consists of those original `resilient' classes that survive, together with others which may arise as a result of the specialization, cf \cite{GerstGiaq:Rigid}.  The family of $q$-Weyl algebras provides an extreme example. 

For brevity, we may refer to a single algebra as being a `deformation' if it is a member of a deformation family and may (despite the foregoing) also refer to a twist as a deformation if it has arisen as in \eqref{f}. 

There is a useful condition insuring that a twist $f$ not be trivial.  If $B$ is a subalgebra of an algebra $A$ then a deformation of $A$ generally does not restrict to one of $B$, as the expression for the deformed product of two elements in $B$ may involve elements of $A$ not contained in $B$. Suppose, however,  that  $\bS'$ is a subsemigroup of  $\bS$, and that $A = \bfk\bS$ has been twisted by $f$.
Since the values of $f$ lie in $\bfk^{\times}$ and do not involve $\bS$, the restriction of $f$ to $\bfk \bS'$  does define a twist of the latter. If the twist of $\bfk \bS$ is trivial, then from the definition one sees that the induced twist of $\bfk \bS'$ must also be trivial. Since it may be easier to prove non-triviality of the latter than of the former, when $A$ is the path algebra of a triangulable space $X$ it will be convenient to choose a $B$ coming from a triangulation of  $X$.

 Henceforth, `cohomology of an algebra', without mention of any module, will mean cohomology with coefficients in itself.
\section{path algebras}\label{path} 
To every small category $\mathcal{C}$ one can associate a semigroup $S\mathcal{C}$ consisting of the morphisms (including identities) of $\mathcal{C}$ and a zero element, the product of two morphisms being their composition when that is possible, and 0 otherwise. Let $\bfk$ be a commutative unital ring. The  semigroup algebra  $\bfk\mathcal{C}$ will be called the \emph{category algebra} or \emph{path algebra} of $\mathcal{C}$ with coefficients in $\bfk$. All concepts of `path algebra' and `path algebra with relations' which occur, e.g. in quiver theory, are special instances.  The objects of $\mathcal{C}$, or more precisely, their identity morphisms, are pairwise orthogonal primitive idempotents.  If $\mathcal{C}$ has only finitely many objects, then the sum of these idempotents is the unit element of $\bfk\mathcal{C}$, which otherwise has no unit element, although there are always ``local units''. That is, for any finitely generated subalgebra $B$ of $\bfk\mathcal{C}$ there is always an idempotent element $e$ such that $ea = ae =a$ for all $a \in B$. Category algebras generalize the concept of poset algebras, which are the category algebras of posets.

In general, there are relations amongst the morphisms in a category. There is a forgetful functor from small categories to quivers which forgets these relations. Its left adjoint is the free category generated by the quiver. The category algebra of a free category is the path algebra of the underlying quiver. There are no relations in this algebra other than that the product of two arrows which can not be concatenated, i.e., such that the starting node of the second differs from the ending node of the first, is zero. In quiver theory it is frequently necessary also to consider ``path algebras with relations'', which essentially just recaptures the concept of a category algebra.  It is simpler to call all of these ``path algebras'' and to call the category algebra of $\mathcal C$ also its path algebra.

A closely related concept is that of the nerve of $\mathcal{C}$, denoted $N\mathcal{C}$, whose $n$-simplices are its $n$-tuples of composable morphisms.  Another restatement of \cite{GS:SC=HC} (cf Theorem \ref{poset}) is the following:
\begin{theorem}
If  the small category $\mathcal{C}$ is a poset, then for any commutative, unital ring $\bfk$ there is a natural isomorphism
$$H^*(N\mathcal{C}, \bfk) \cong H^*(\bfk\mathcal{C}, \bfk\mathcal{C}),$$
where on the left one has simplicial cohomology and on the right Hochschild cohomology.  $\square$
\end{theorem}
\noindent Since the barycentric subdivision of a simplicial complex is a  poset and also has the same cohomology as the original complex, one has the theorem stated in the title of \cite{GS:SC=HC}. For a general small category $\mathcal{C}$, the relation between $H^*(N\mathcal{C}, \bfk)$ and $H^*(\bfk\mathcal{C}, \bfk\mathcal{C})$ is given by the Cohomology Comparison Theorem, \cite{GS:Monster}.

The usual definition of a path in a topological space $X$ is a continuous map $\g: [0,1] \to X$, the product of two such,  $\g$ and $\g'$, being  defined as follows.  If $\g(1) \ne \g'(0)$ then set $\g\g'=0$.  Otherwise set 
\begin{equation*}
(\g\g')(t) =
		\begin{cases}
				\g(2t)			& \mbox{ if $0 \le t\le 1/2$} \\
				\g'(2(t-1/2))	& \mbox{ if $1/2 \le t \le 1$}.
		\end{cases}
\end{equation*}
With this definition,  however,  the product is not associative; the two ways of associating a concatenation of three paths, namely, $(\g\g')\g''$ and $\g(\g'\g'')$ are not identical but only homotopic. 
(That a specific homotopy would have non-trivial implications for relating the five ways of associating four paths  was one of the observations that led Stasheff  to the definition of $A_{\infty}$ and $L_{\infty}$ algebras.) 
To make paths into a strictly associative semigroup, following  \cite{Moore:Paths}
define a \emph{Moore path} on a topological space $X$ to be a pair $(\g, \ell)$, where $\ell \in [0,\infty)$ and $\g:[0,\infty) \to X$ is a continuous mapping such that $\g(t) = \g(\ell)$ for $t \ge \ell$. The \emph{start} of $\g$ is $\g(0)$,  its \emph{length} is $\ell$, and its \emph{end} is $\g(\ell)$.   When  the start of a path $(\g', \ell')$ coincides with the end of $(\g, \ell)$ define their concatenation, $(\g\g', \ell+\ell')$, by setting 
\begin{equation*}
(\g\g')(t) =
		\begin{cases}
				\g(t)		& \mbox{ if $0 \le t\le \ell$} \\
				\g'(t-\ell)	& \mbox{ if $\ell \le t < \infty$}
		\end{cases}
\end{equation*}
This product is associative and is continuous on the set of pairs of paths that can be concatenated. Defining the  product of paths which can not be concatenated to be zero, one has a semigroup with zero, $\bS X$.   The \emph{path algebra $\bfk X$ of $X$}  is the semigroup algebra generated over the commutative unital ring $\bfk$ by this semigroup. The set $X^Y$ of all mappings from a space $Y$ to a space $X$ carries the natural compact-open topology, so there is a natural topology on the space $X^{[0,\infty)}\times [0,\infty)$ of paths. 
The path space of $X$ can be continuously retracted to $X$, so its homotopy and singular homology groups are identical with those of $X$ itself. 

On a Riemannian manifold $\cM$, `path'  will always mean a piecewise smooth mapping. Such paths can be parameterized by arc length; any path which is so parameterized and which has Riemannian arc length $\ell$ will be considered as a Moore path of length $\ell$. The concepts of Moore path length and Riemannian arc length will then coincide.  There is  a standard metric on the space of piecewise smooth maps $\g: [0,1] \to \cM$, cf e.g.,  
\cite[p. 88]{Milnor:Morse}.  
This can be adapted to piecewise smooth Moore paths $(\g_1, \ell_1), (\g_2, \ell_2)$ by defining the distance between them to be
\begin{multline}\label{Moore metric}
d((\g_1, \ell_1), (\g_2, \ell_2)) = \\ \max_{0\le t < \infty}d(\g_1(t), \g_2(t)) + \text {\Large[}\int_0^{\infty}(\frac{ds_1}{dt} -\frac{ds_2}{dt})^2dt \,\text {\Large]}^{1/2}  +|\ell_2-\ell_1|;
\end{multline}
both limits at $\infty$ can be replaced by $\max(\ell_1, \ell_2)$. 

 Henceforth we will tacitly assume that any Riemannian manifold $\cM$ is connected and complete, and unless otherwise specified, that a piecewise geodesic path on $\cM$ is parameterized by its arc length. Mention of its length can then be omitted. Those paths which are piecewise geodesic form a subalgebra of the full path algebra.

For an arbitrary topological space $X$, multiplication in the semigroup $\bS X$  is generally not continuous, since a sequence of pairs of non-concatenatable paths, whose products are therefore always 0, can converge to a pair which can be concatenated.  However,
on a Riemannian manifold $\cM$  one can define (but we will not need) a continuous product in the full semigroup algebra $\R\cM$: Choose a smooth, monotone decreasing function $f:[0,1] \to \R$ with $f(0) = 1, f(1) = 0$.  Suppose that $(\g, \ell)$ is a path with end $q$ and that $(\g', \ell')$ is one whose start is $q'$.  If the radius of injectivity of $q$ is $\rho$ and  $d(q,q') \ge\rho$, then define their product to be 0.  Otherwise, let $\g_{q,q'}$ be the unique minimal geodesic which starts at $q$ and ends at $q'$. If its length is $\ell''$, redefine $\g\g'$ to be $f(\ell''/\rho)(\g\g_{q,q'}\g', \ell+\ell'+\ell'')$. 

\section{Taille}\label{taille}  The $k$th homotopy group of a topological space $X$ with base point, $\pi_k(X)$, is the group of  homotopy classes of base point preserving mappings of the  $k$-sphere $S^k$ into $X$.  For $k =1$ it is the classical, generally non-Abelian, fundamental group but for $k  \ge 2$ it is always Abelian. We will be most interested in the case $k = 2$.   One can have $\pi_k(X) = 0$ 
while the $k$th integral homology group  $H_k(X)$ of $X$
does not vanish. This is the case, for example, with the 2-torus $T^2$ since $\pi_2(T^2) = 0$ but $H_2(T^2) = \Z$.  By contrast, $H_3(S^2) = 0$, while $\pi_3(S^2) = \Z$. The latter group  is  generated by the Hopf fibration of $S^3$ over $S^2$. 

The Hurewicz theorems provide a key link between homotopy groups and homology groups, but of them we need here only  that on any space $X$ one has a morphism $\pi_k(X) \to H_k(X)$. Choosing  any Abelian coefficient group $G$,  the Universal Coefficient Theorem provides an epimorphism $H^k(X,G) \to \operatorname{Hom}(H_k(X), G)$. With the Hurewicz theorems,
 it follows that every $k$-cocycle $\omega$ of $X$ with coefficients in $G$  defines a morphism $\pi_k \to G$ which depends only on the class of $\omega$  in $H^k(X, G)$. The  \emph{taille group} of $\omega$ is defined to be the image group, denoted $\omega(\pi_k(X))$. If $G = \Z$, then $\omega(\pi_k(X))$ must be  cyclic and we define the \emph{taille} of $\omega$ to be zero if $\omega(\pi_k(X))$ is trivial and otherwise to be its positive generator.  In particular,  a closed $k$-form $\omega$  on a differentiable manifold  $\cM$ gives rise to a morphism of $\pi_k(\cM)$ into $\R$.

\begin{defn} Let $\omega$ be a closed $k$-form on a differentiable manifold $\cM$.  The \emph{taille} of $\omega$ is i) 0 if $\omega(\pi_k(\cM)) = 0$,  ii) the positive generator of $\omega(\pi_k(\cM))$ if the latter is a non-trivial cyclic subgroup of $\R$, or iii) $\infty$ otherwise. If $\cM$ is symplectic then by its \emph{taille}, $\tau(\cM)$, we will mean the taille of its defining symplectic 2-form.
\end{defn} 
 The exterior powers of the defining form of a symplectic manifold map the  corresponding even dimensional homotopy groups into $\R$, so there are secondary tailles for each even dimension greater than two.  

A closed path or loop $\g$ in $\cM$ is contractible if there is a continuous map  of a 2-disc into $\cM$ the restriction of whose image to the boundary is the image of the path; we will call this a spanning disk.  A closed 2-form $\omega$ can be integrated over the image of the disc, but the result is not a function of $\g$ alone but also of the choice of spanning disk.  If we have a second spanning disk then, together with the first, one has  a map of the 2-sphere into $\cM$. The difference between the values obtained will be the integral of $\omega$ over the image of the 2-sphere. If $\omega$ has a finite taille, $\tau$, then reduction modulo $\tau$ gives a well-defined value in $\R/\tau$, independent of the choice of spanning disk. This will be denoted by $\langle \omega, \g \rangle$.  If $\cM$ has a finite fundamental group,
 then for any loop $\g$ there is some positive multiple $m\g$ which is contractible
and we can set $\langle \omega, \g \rangle = \langle \omega, m\g \rangle/m$.  

Let $\cM$ now be a Riemannian manifold and $\omega$ be a closed 2-form on it with  finite taille $\tau$. We want to define an additive coherent 2-cocycle $\tilde\omega$ of the path semigroup  $\bS\cM$ of $\cM$ whose value $\tilde\omega(\gamma,\gamma')$ will depend only on the homotopy classes of its arguments $\gamma,\gamma' $.  If  $\gamma$ and  $\gamma'$ can not be concatenated,
 then set $\tilde\omega(\gamma,\gamma') = 0$. If they can, then $\tilde\omega(\gamma,\gamma')$ will be left undefined unless i) there is a unique shortest geodesic $\z$ homotopic to $\gamma$, ii) there is a unique shortest geodesic $\z'$ homotopic to $\gamma'$ and iii) there is a unique shortest geodesic $\z''$ homotopic to their concatenation $\gamma\gamma'$. When the conditions are met,  $\z''$ must also be homotopic to the concatenation of $\z$ and $\z'$, and the homotopy then provides an element of area bounded by $\z, \z'$ and $\z''$, oriented by following 
them in that order. Set $\tilde\omega(\gamma,\gamma')=\int \omega$ over this element, taken modulo the taille, $\tau$ (which may be zero).
Note that the integral defining $\tilde\omega(\gamma,\gamma')$ is over a triangle whose sides are uniquely defined geodesics which depend only on the homotopy classes of $\gamma$ and $\gamma'$, but the homotopy defining the element of area is not unique. If we have two distinct homotopies then, as before, they  define a mapping of the 2-sphere $S^2$ into $\cM$. The difference between the integrals will be zero if this mapping is homotopic to zero but possibly not otherwise. However,  the difference between the integrals will be a multiple of $\tau$, and the integral  becomes well-defined and independent of the choice of homotopy if reduced modulo $\tau$. (For the related concept of a homotopy multiplicative map, 
cf, eg, \cite{Stasheff:Homotopy associative}.)

That $\tilde\omega(\gamma,\gamma')$ is sometimes undefined suggests  enlarging the concept of a topological algebra to allow products to be undefined when there is some measure of smallness for the set of cases where this occurs, and similarly for morphisms, cochains, and other constructs. 
Since $\cM$ carries a volume form, and hence a measure, so does $\cM \times \cM$.  It seems reasonable to require that the set of pairs of points $a,b \in \cM$ such that some homotopy class of paths between them contains no unique shortest geodesic has measure zero and does not disconnect $\cM \times \cM$, something which we will henceforth tacitly assume.

\begin{defn}
The \emph{tailleur cochain} of a closed de Rham 2-form  $\omega$ on a Riemanian manifold $\cM$ is its associated additive 2-cochain $\tilde\omega$ of the path algebra of $\cM$.
\end{defn}

\begin{theorem} The tailleur cochain of $\cM$ relative to a closed de Rham 2-form $\omega$ is a cocycle. If the taille, $\tau$, of $\cM$ relative to $\omega$ is zero, then the tailleur cocycle $\tilde\omega$  defines a one-parameter family of coherent deformations of the real path algebra $\R\cM$ with multiplicative twisting 2-cocycle $f_{\hbar} = \exp(\hbar\tilde\omega)$ whose values are always strictly positive. If $\tau > 0$, then the tailleur cocycle defines a single coherent twist of the complex path algebra $\C\cM$ with twisting cocycle  $f = \exp((2\pi i/\tau)\tilde\omega)$ whose values lie on the unit circle. 
\end{theorem}
\noindent\textsc{Proof.} Suppose that we have three paths $\gamma, \gamma', \gamma''$ which can be concatenated to give $\gamma \gamma'\gamma''$, and consider $\tilde\omega$ for the moment as having values, as originally, in $\R$ rather than in $\R/\tau\R$. The coboundary of $\tilde\omega$ evaluated on these three paths is 
\begin{equation*}
\delta \tilde\omega(\gamma,\gamma',\gamma'')=\tilde\omega(\gamma',\gamma'') -\tilde\omega(\gamma\gamma',\gamma'') +\tilde\omega(\gamma,\gamma'\gamma'')  -\tilde\omega(\gamma,\gamma'),
\end{equation*}
cf \eqref{coboundary}. This is just the integral of the closed 2-form $\omega$ over the surface of a 3-simplex in $\cM$ (which happens to have geodesic edges).   It is therefore a multiple of the taille, so  $\tilde\omega$  becomes an additive 2-cocycle after reduction modulo  $\tau$.  The rest follows. $\square$
\smallskip

When $\tau = 0$ the tailleur cocycle $\tilde\omega$, whose values initially lie in $\R$, can be reduced modulo any $\tau >0$. With coefficients extended to $\C$, the reduced cocycle will produce
a single coherent twist, but as observed previously, the deformations induced will be entirely different.


The definitions above apply also to the \emph{geodesic path algebra} of $\cM$, namely, the subalgebra generated by all piecewise geodesic paths on $\cM$. The \emph{homotopy path algebra} of $\cM$ is generated by triples $(a,b, [\gamma])$, where $(a,b)$ is an ordered pair of not necessarily distinct points of $\cM$ and $[\gamma]$ is a homotopy class of paths $\gamma$ from $a$ to $b$.  The product $(a,b,[\gamma])(b',c,[\gamma'])$ is zero if $b \ne b'$, and otherwise is $(a,c,[\gamma\gamma'])$. The path algebra maps onto the homotopy path algebra by sending every path to its homotopy class. Since the value of $\tilde\omega$ on a pair of concatenatable paths $\g, \g'$ depends only on their homotopy classes, the deformation that the tailleur cocycle induces on the path algebra descends to the homotopy path algebra. The geodesic path algebra is `almost' isomorphic to the  homotopy path algebra,  since it is just the homotopy path algebra with the multiplication left undefined when there is no unique shortest geodesic in $[\gamma\gamma']$, so a deformation is also induced on it. The definition of the homotopy path algebra is, like that of the path algebra, independent of the metric on $\cM$, but the set where the almost isomorphism is undefined may depend on the metric, as may any induced twist.  

 Every differentiable manifold $\cM$ has a triangulation, J. H. C. Whitehead \cite{JHCWhite:Triang},  i.e, there exists  a simplicial complex $K$ together with a homeomorphism $\theta:K \to \cM$. The cohomology of $\cM$ can then be identified with the simplicial cohomology of  $K$. This is based on earlier work of S.S. Cairns; for a brief history cf Saucan \cite{Saucan:Munkres}. A non-differentiable manifold need not have a triangulation but  piecewise linear manifolds, although generally not smoothable, of course do. 
 
 The barycentric subdivision  $K'$ of $K$ has the same cohomology as $K$ and is also a poset whose  1-simplices can be parameterized so that as Moore paths, each has length 1. Its image under $\theta$ is therefore a poset subalgebra, denoted $ A(\theta)$, of the path algebra $k\cM$ of $\cM$. 
The de Rham cohomology of $\cM$ is then naturally isomorphic to $H^*(A(\theta),A(\theta))$.  In view of the remark at the end of \S2 we therefore have the following.

\begin{theorem}\label{nontrivial} Let $\cM$ be a Riemannian manifold and $\omega$ be a closed 2-form with finite taille $\tau$.  Then the coherent twist of the path algebra induced by $\omega$ is trivial if and only if the class of $\omega$ is trivial as an element of $H^2(\cM, \mathbb R/\tau\R)$. $\square$
 \end{theorem}
 \vspace{-2mm}

\section{Quantization of Phase Space}\label{dB}

A phase space whose underlying configuration space is Euclidean simultaneously carries Euclidean, orthogonal, and symplectic structures in which the spatial coordinates are orthogonal to the momentum coordinates, cf Arnol'd, \cite{Arnol'd:Maslov}. The relationship between these structures allowed Maslov to assign an integer-valued index to every smooth closed curve in such a space.  Arnol'd gives an elegant redefinition of this Maslov index under which the index of a smooth simple closed curve in the phase space of a particle moving on a line is easily seen to be 2 because its tangent line, which is considered undirected, turns twice, cf \cite{Arnol'd:Maslov}.

Mathematically, phase space is the cotangent bundle of configuration space, so its homotopy groups are just those of the latter. They may not vanish, but the generators of its $\pi_2$ are all completely isotropic relative to the canonical symplectic form $\omega$. Therefore, if $\xi \in \pi_2$ then $\omega(\xi)= 0$, so the taille of mathematical phase space is always zero. Nevertheless, the existence of a quantum of action, wave-particle duality, and Bohr's  model of the hydrogen atom (even if not completely successful),  strongly suggest that the phase space in which a particle actually moves has positive taille equal to Planck's original constant, $h$. As observed in \S \ref{twist}, an additive 2-cocycle with coefficients in $\R$ also defines one with coefficients in $R/\tau$ for any $\tau > 0$, but the corresponding twist is quite different. The taille of mathematical phase space is zero, so its canonical additive tailleur  2-cocycle can be reduced modulo any $\tau > 0$ to give a quantization of its path algebra. For the phase space of a particle moving on a straight line, setting $\tau$ equal to $h$  implies, as we will show, wave-particle duality.

The phase space of a particle moving on a straight line has momentum and position coordinates $p$ and $q$, and canonical symplectic form $\omega = dp\wedge dq$.   As the Maslov index of a non-singular simple closed curve in this space is 2, we choose $2\omega$ (rather than just $\omega$) to define a coherent twist of the path algebra,  \cite{Dolgushev:private}. Suppose now that a particle at rest at the origin is given an impulse at time 0 which instantaneously gives it momentum $p$, as a result of which it then moves freely with constant velocity $v$. (This can also be read in reverse, as when a particle is stopped.) Its path in phase space up to time $t$ is then a broken line consisting of a segment $\g$ of length $p$ in the momentum direction, followed by an orthogonal segment $\g'$ of length $vt$ in the spatial direction ending at $(p, vt)$. To compute the values of the tailleur cocycle $\tilde\omega(\g,\g')$, following the procedure of the preceding section,
 we must 1) complete this path to a closed curve returning to the starting point using a geodesic homotopic to the concatenation $\g\g'$ and 2) integrate $\omega$ over the spanning disk defined by the homotopy.  Here the straight line $\g''$ back to the start  is clearly homotopic to $\g\g'$, the spanning disk is the right triangle whose sides are $\g,\g',\g''$, and the integral in question is its area, $(1/2)pvt$.  This area is an action integral, but it must be taken modulo the taille and multiplied by the Maslov index of the closed curve $\g\g'\g''$.

One can not immediately assign Maslov index 2 to this curve because of its singularities at the corners. It must be approximated by a smooth curve, and as the curve varies through the space  of paths, it must be that if the approximation is sufficiently close in the given topology, then the index of the smoothed curve will be independent of the smoothing.  Closeness in the metric on Moore paths introduced in \eqref{Moore metric} is alone not sufficient to insure that smooth curves near a given one have the same Maslov index. In phase space, however, flow lines never cross.  Adding the condition that approximating curves do not introduce additional crossings, any smooth curve close  to the original triangle will have Maslov index 2.

With this, the value of the  additive tailleur 2-cocycle $\tilde\omega$ on the pair of concatenatable paths $\g, \g'$  is just twice the area of the right triangle they define, i.e., $pvt$,  so the value of the multiplicative twisting 2-cocycle $f(\g,\g')$ is $\exp({(2\pi i/h)pvt)}$. The deformed  product of the two segments $\g, \g'$ is therefore
$$\g\star\g' = \exp({(2\pi i/h)pvt)}\cdot\g \g'  = \exp({(i/\hbar)pvt)}\cdot\g \g'. $$
Here $\g$ is fixed, $\g'$ depends on $t$, and the twisting function $\exp({(2\pi i/h)pvt)})$ is periodic in time with period $h/pv$, in which time the particle moves a distance $v(h/pv)$. Its wave length is therefore that hypothesized by de Broglie's wave-particle duality, namely $h/p$.  Wave-particle duality thus appears as a deformation of the path algebra of phase space. One must integrate over an increasingly large area as the particle travels, indicating that the wave associated to a particle by this argument is not local since it does not depend only on the nature of phase space in the immediate vicinity of the path of the particle. Because of the flatness of the present phase space, which is just a plane, the wave it associates to the moving particle has a wavelength which is constant in time. 

The preceding argument required  reduction of the values of the tailleur cocycle modulo $h$, which suggests that phase space must somehow already have a positive taille, i.e., must be quantized.  Mathematical phase space, however, always has taille equal to zero, the basic paradox mentioned in the introduction.    One resolution is to look for a larger `physical' phase space in which the particle actually moves,  one having taille equal to $h$, and of which mathematical phase space is a projection in order to pull back the canonical symplectic form. However,  this larger physical phase space must then have singularities.  

For a particle  moving on a line, a simple construction of a space with the necessary properties would be to take two copies of mathematical phase space, a plane, tile both the same way with tiles of area $h/2$, and identify the tile edges in one copy with the corresponding edges in the other.  The most easily visualized tiling is rectangular with edges aligned with the space and momentum directions, but the physical phase space in which a particular particle moves might depend on its charge and velocity and might not even be flat. The nature of the tiling is so far unknown, but whatever it may be, this physical phase space does have a projection onto mathematical phase space which is 2 to 1 except where the edges of the tiles have been identified (thereby producing singularities), where it is 1 to 1. Elsewhere, it acquires a defining 2-form by pull back of that of the mathematical phase space and this 2-form can be integrated over elements of area because the singularities have measure zero when projected back to mathematical phase space. In view of the double covering, one should  count areas in mathematical phase space as doubled, something which the Maslov index does by viewing tangent lines, even to an oriented curve, as unoriented. This space has taille  equal to $h$, as each doubled tile, which is homeomorphic to a 2-sphere, has total area $h$. 

The space that has been constructed  is homeomorphic to a lattice of spheres which are disjoint except that the equator of each is divided into several parts, with each part identified with a part of the equator of another; a close representation  would be a sheet of bubble wrap. Except at the singularities the space has a `top' and a `bottom', interchange of which  reverses orientation.  This introduces spin: A point on top is in `spin up' position, a point on the bottom has `spin down', while at the singular points spin is undefined. Distinct particles can have identical positions and momenta as long as they have different spins. (A 2-sphere is homeomorphic to complex 1-dimensional projective space, $\C P^1$; when the configuration space of the particle is $n$-dimensional, one should probably replace $\C P^1$ by $\C P^n$, since $\pi_2(\C P^n) = \Z$  for all $n \ge 1$.)

This model of physical phase space in effect gives phase space itself a dual nature.  At the macroscopic level it is indistinguishable from its projection onto mathematical phase space; it is just a symplectic manifold.  At the Planck level, however, it is a lattice of small tiles or cells, and any physical particle should probably be viewed as large compared to a single cell.  While the nature of the tiling remains  unknown, it seems reasonable to conjecture that at any given time a particle completely occupies either a single cell or some fixed number of contiguous cells.  A photon or neutrino might occupy but a single cell, while a massive  particle will occupy a number of contiguous cells proportional to its mass,  the proportionality factor being dependent on the still uncertain nature of the tiling. The spread of a packet in time can then be described as a random walk. In finite time there could therefore be only be a finite number of Feynman paths. Consider first the case of a moving particle occupying at any give time but a single cell and moving with preservation of spin (remaining always either on the top tile or the bottom one, unless acted upon externally). Suppose that at each unit of time it has probability $\wp$ of moving one unit to the right and $1-\wp$ of moving one unit to the left. Its average rightward velocity or drift is then $2\wp - 1$.  Let  $p$ now denote the momentum of the particle, averaged over the single cell it occupies, and let $m$ be its mass. Then $v = p/m$ is its mean velocity.  Set $\overline\upsilon = p/mc$, where $c$ is the speed of light; this is the drift normalized with speed of light set equal to 1, a number between $-1$ and $+1$.  Conversely, if the drift $\overline\upsilon$ is given then  $\wp = (\overline\upsilon +1)/2$.
 
Choosing as the unit of time the average time it takes for light to cross a cell (which will depend on the tiling),
in this model a particle can never travel faster than the speed of light, while in the Schr{\"o}dinger picture a particle would always have a very small but positive probability of being farther from its start at any positive time than the speed of light would allow. In this respect the Schr{\"o}dinger model resembles the Gaussian approximation to the binomial distribution. The Gaussian distribution has infinitely long tails but is a reasonably accurate approximation to the binomial distribution as long as one stays within two standard deviations of the true end points of the binomial distribution. The relation of the Schr{\"o}dinger model  to the discrete model presented here may be similar to that of the Gaussian distribution to the binomial.  For a particle moving at the speed of light we have $\wp = 1$ and $1-\wp =0$, so for such a particle the distribution does not spread at all.  For a massive particle occupying $N$ contiguous cells, each cell must have the same drift, giving 
 $N$ identical but independent probability distributions, of which we must take the mean. For a particle with mean momentum zero, the distribution of the particle's possible locations after a large number of units of time should, with appropriate choice of tiling, be asymptotically the same as that given by the Schr{\"o}dinger picture.  

In $\S \ref{dB}$, phase space was assumed to be flat, and as a result, the wave associated to a particle had constant frequency.  As an example of what might happen to a particle moving in a curved phase space, consider the deformation of the path algebra of the unit sphere $S^2$ induced by the usual area form. While it is assumed to be tiled at the Planck level, macroscopically we may still view it as a manifold. Geodesics between antipodal points are not unique, so the product of concatenatable paths will be undefined when the beginning of the first is antipodal to the end of the second. However, the set of pairs of antipodal points in $S^2 \times S^2$ has codimension two and does not disconnect that space. The usual area form $\omega$ is, up to constant multiple, the only harmonic $2$-form and its class spans the second de Rham cohomology module of the sphere. The taille of $\omega$ is the area of the sphere, $4\pi$, and the twist it defines is non-trivial by Theorem \ref{nontrivial}.  

Imagine now that this were the physical phase space seen by some kind of particle, and view the meridians of longitude as lines along which the spatial dimension is constant and momentum varies.  Suppose that a particle in the southern hemisphere at 0 degrees longitude receives an impulse which instantaneously translates it north on that meridian to the equator, and thereafter moves eastward.  From any point on the particle's equatorial motion after the initial impulse, complete its path to a closed curve by the shortest geodesic to the point of origin.  As the particle continues to circle the equator, it will experience a change in phase which is no longer linear in time since the area of the spherical triangle swept out will not be linear in time. The particle's wavelength will therefore also change, although its velocity remains constant and there is no change in the local geometry. That the observed frequency can change with no change in energy suggests that part of the energy of a particle moving in a curved phase space might not be directly visible.
 
In a curved physical phase space, the path of a particle moving freely after an initial impulse can also pass through points from which there is no unique shortest geodesic to the point from which it started. These will generally be points of discontinuity of its wave function, and therefore ones where the particle might appear to undergo a discontinuous change. 

The quantized model of physical phase space constructed here seems to be the simplest one resolving the basic paradox that mathematical phase space presents, but is surely not the only one.  In others, possibly projecting onto the present one,  the points might represent further measurable physical attributes of a particle.   It seems certain, however,  that any model which gives physical phase space an intrinsic taille equal to $h$, something  necessary for agreement  with de Broglie's wave-particle duality, will not only give rise to spin, but will have further physical implications.


\begin{thebibliography}{10}

\bibitem{Arnol'd:Maslov}
V.~I.~Arnol'd. 
\newblock{On a characteristic class entering into conditions of quantisation}. (Russian. English translation) 
\newblock{\em Functional Analysis and Its Applications}, 1:1--14, 1967.

\bibitem{FiniteRep}
R.~Bautista, P.~Gabriel, A.~V.~Roiter, and L.~Salmeron.
\newblock {Representation-finite algebras and multiplicative bases}.
\newblock {\em Inventiones Math.}, 81:217--285, 1985.

\bibitem{BFFLS}
F.~Bayen, M.~Flato, C.~Fr{\o}nsdal, A.~Lichnerowicz, and D.~Sternheimer.
\newblock {Deformation theory and quantization, I and II.}
\newblock{\em Ann. of Phys.}, 111:61--151, 1978.

\bibitem{DitoSternheimer:Genesis}
G.~Dito and D.~Sternheimer.
\newblock{Deformation quantization: genesis, developments
and metamorphoses}, pp. 9--54 in
\newblock{G. Halbout, ed., \textsc{Deformation Quantization, IRMA Lectures in Math. Theoret. Phys.
1},}
\newblock{Walter de Gruyter, Berlin} 2002.
\newblock{\texttt {arXiv:0201168v1 [math.QA] 18 Jan 2002}}.

\bibitem{Dolgushev:private}
V.~Dolgushev.
\newblock{Private communication.}

\bibitem{G:Cohomology}
M.~Gerstenhaber.
\newblock{The cohomology structure of an associative ring}.
\newblock{\em Ann. of Math.} 78:267--288, 1963.

\bibitem{G:Def1}
M.~Gerstenhaber.
\newblock {On the deformation of rings and algebras}.
\newblock {\em Ann. of Math.}, 79:59--103, 1964.

\bibitem{G:Def3}
M.~Gerstenhaber.
\newblock {On the deformation of rings and algebras:III}.
\newblock {\em Ann. of Math.}, 88:1--34, 1968.

\bibitem{G:Least uncertainty}
M.~Gerstenhaber.
\newblock{Least uncertainty principle in deformation quantization}. 
\newblock{\em J. Math. Phys.} 48:022103,15 pp., 2007. 

\bibitem{G:Self-dual}
M.~Gerstenhaber.
\newblock{ Self-dual and quasi self-dual algebras}.
\newblock{\em  Israel J. Math.}, 200:193--211.
\newblock {\texttt{arXiv:1111.0622v1 [math.KT] 2 Nov 2011}}.


\bibitem{GerstGiaq:Rigid}
M.~Gerstenhaber and A.~Giaquinto.
\newblock{Deformations associated with rigid algebras}.
\newblock{\em J. Homotopy Relat. Structures}.
\newblock{DOI 10.1007/s40062-013-0068-x}.
\newblock{{\texttt{arXiv:1208.5068v1 [math.QA] 24 Aug 2012}}}.


\bibitem{GS:SC=HC}
M.~Gerstenhaber and S.~D. Schack.
\newblock{Simplicial cohomology is Hochschild cohomology}.
\newblock {\em J. Pure and Appl. Algebra}, 30:143--156, 1983.

\bibitem{GS:Hodge}
M.~Gerstenhaber and S.~D.~Schack.
\newblock{A Hodge-type decomposition for commutative algebra cohomology}.
\newblock {\em J. Pure and Appl. Algebra}, 48:229--247, 1987.

\bibitem{GS:Monster}
M.~Gerstenhaber and S.~D. Schack.
\newblock {Algebraic cohomology and deformation theory}, in 
 {\textsc {Deformation Theory of Algebras and Structures and
  Applications}} (M.~Hazewinkel and M.~Gerstenhaber, eds.), pp. 11--264.
 \newblock Volume 247 of NATO ASI Science Series.
\newblock Kluwer Academic Publishers, Dordrecht/Boston/London, 1988.

\bibitem{GiaquintoZhang:Bialgebra}
A.~Giaquinto and J.~J.~Zhang.
\newblock{Bialgebra actions, twists, and universal deformation formulas}.
\newblock {\em J. Pure and Appl. Algebra}, 1288:133--151, 1998.

\bibitem{Groenewold}
H.~J. Groenewold.
\newblock {On the principles of elementary quantum mechanics}.
\newblock {\em Physica}, 12:405--460, 1946.

\bibitem{HKR}
G.~Hochschild, B.~Kostant, and A.~Rosenberg.
\newblock{Differential forms on regular affine algebras}.
\newblock{\em Trans. Amer. Math. Soc.}, 102:383--408, 1962.

\bibitem{Kontsevich:Poisson}
M.~Kontsevich.
\newblock{Deformation quantization of Poisson manifolds}.
\newblock{\em Lett. Math. Phys.} 66:157--216, 2003.

\bibitem{Loday:Euler}
J.-L.~Loday.
\newblock{Partition eul{\'e}rienne et op{\'e}rations en homologie cyclique}. (French. English summary)
\newblock{\em C. R. Acad. Sci. Paris S{\'e}r. I Math.} 307:283--286, 1988.

\bibitem{Milnor:Morse}
J.~Milnor.
\textsc{Morse Theory} (Based on lecture notes by M.~Spivak and R.~Wells).
\newblock{Annals of Mathematics Studies Number 51}.
\newblock Princeton University Press, Princeton, New Jersey, 1963.

\bibitem{Moore:Paths}
J.~C.~Moore.
\newblock{ Le th\'eor\`eme de {F}reudenthal, la suite exacte de {J}ames et l'invariant de {H}opf g\'en\'eralis\'e}, in \textsc{  S\'eminaire Henri Cartan 7 no. 2, 1954-55}, Expos{\'e} 22, 15 pp.
\newblock{C. B. R. M., Gauthiers-Villars}, Paris 1955.

\bibitem{Moyal}
J.~E.~Moyal.
\newblock {Quantum mechanics as a statistical theory}.
\newblock {\em Proc. Cambridge Phil. Soc.}, 45:99--124, 1949.

\bibitem{Novikov:Semigroup}
B.~V.~Novikov.
\newblock {Semigroup cohomology and applications}.
\newblock{em Algebra -- Representation Theory}, (2001): 219--234.
\newblock{\texttt{arXiv:0803.0463v1 [math.RA] 4 Mar 2008}}.

\bibitem{Rieffel:actions}
M.~A.~Rieffel.
\newblock{Deformation Quantization for Actions of $\bR^d$}.
\newblock{\em Memoirs Amer. Math. Soc.} vol. 106, no. 506, November, 1993.

\bibitem{Romero-Ayala}
R.~ Romero-Ayala,  C.~ A.~ Soto-Campos, and S.~ \mbox{Valdez-Alvarado}.
\newblock{A noncommutative embedding of Reissner-Nordstr{\o}m Spacetime}.
\newblock{\texttt{arXiv:1503.05469v1 [gr-qc] 18 Mar 2015}}.

\bibitem{Saucan:Munkres} 
E.~Saucan.
\newblock {Note on a theorem of Munkres}.
\newblock{\em Mediterranean Journal of Mathematics}, 2:215--229, 2005.  
\newblock {\texttt{arXiv:math/0403055v2 [math.GT] 14 Mar 2004}}.

\bibitem{Stasheff:Parallel}
J.~D.~Stasheff.
\newblock{``Parallel'' transport - revisited}.
\newblock \texttt{arXiv:11115495v1 [math.AT] 23 Nov 2011}.

\bibitem{Stasheff:Homotopy associative}
J.~D.~Stasheff.
\newblock{Homotopy associativity of H-spaces, I, II}.
\newblock {\em{Trans. Amer. Math. Soc.}}, 108:275-292, 294--312, 1963.

\bibitem{Sternheimer:20YearsAfter}
D.~Sternheimer.
\newblock{Deformation quantization: Twenty years after}, in \textsc{Particles, Fields, and Gravitation}, {\L}odz, 1998 (J. Rembieli{\'n}ski, ed.).
\newblock AIP Press, New York, 1998, 107--145.
\newblock \texttt{arXiv:math/9809056 [math.QA] 10 Sep 1998}.

\bibitem{Waldmann:Poisson-Geometrie}
S.~Waldmann.
\newblock{\textsc{Poisson-Geometrie und Deformationsquantisierung}}.
\newblock{Springer Berlin Heidelberg New York}, 2007.


\bibitem{JHCWhite:Triang}
J.~H.~C.~Whitehead.
\newblock {\em On {$C^1$-complexes}}. 
\newblock {\em{Ann. of Math.}}, 41:809–-824, 1940.

\end{thebibliography}
\end{document}